\documentclass[10pt]{article}
\usepackage{amsfonts}
\usepackage{amsmath}
\usepackage{amssymb}
\usepackage{amsthm}
\usepackage{algorithm}
\usepackage{algorithmicx}
\usepackage{algpseudocode}
\usepackage{graphicx}
\usepackage{subcaption}
\usepackage{indentfirst}
\usepackage{placeins}
\usepackage{flafter}
\usepackage{cite}
\usepackage[hidelinks]{hyperref}
\usepackage{xurl}

\numberwithin{equation}{section}
\theoremstyle{plain}
\newtheorem{theorem}{Theorem}[section]
\newtheorem{lemma}[theorem]{Lemma}
\newtheorem{proposition}[theorem]{Proposition}
\newtheorem{corollary}[theorem]{Corollary}
\theoremstyle{remark}
\newtheorem{remark}[theorem]{Remark}
\numberwithin{algorithm}{section}
\algrenewcommand\algorithmicrequire{\textbf{Input:}}
\algrenewcommand\algorithmicensure{\textbf{Output:}}

\newcommand{\R}{\mathbb{R}}

\newcommand{\eps}{\varepsilon}
\newcommand{\rank}{\operatorname{rank}}
\newcommand{\range}{\operatorname{range}}

\begin{document}
\hypersetup{
 pdftitle={Robust High-Order BUG Integrators without Basis Augmentation},
 pdfauthor={Shiheng Zhang}
}
\title{Robust High-Order BUG Integrators\\
without Basis Augmentation}
\author{Shiheng Zhang\thanks{Department of Applied Mathematics, University of Washington, Seattle, WA 98195, USA (\texttt{shzhang3@uw.edu}).}}
\date{}
\maketitle

\begin{abstract}
We develop robust high-order basis-update and Galerkin (BUG) integrators
for dynamical low-rank approximation without basis augmentation or rank
truncation. The basis and core updates use the same approximate increment.
We establish a two-sided relaxed-exactness estimate that controls the
error in the low-rank approximation by the increment error. Under the
stated assumptions on the vector field and its normal component,
sufficiently accurate increments yield local and global error bounds
whose constants and step-size restrictions are independent of small
singular values. Using explicit Runge--Kutta increments, we obtain robust orders two,
three, and four for the explicit midpoint method, Heun's third-order
method, and Merson's fourth-order method, respectively. Every internal
stage and final approximation has rank at most the prescribed rank $r$.
We also prove that these orders are retained under arbitrary basis
augmentation followed by best rank-$r$ SVD truncation at the internal
stages and the endpoint, without a separate assumption on the truncation
error. Numerical experiments confirm the convergence orders
and show comparable performance with and without augmentation.
\end{abstract}

\medskip
\noindent\textbf{Keywords.}
Dynamical low-rank approximation; basis-update and Galerkin integrator;
relaxed exactness; Runge--Kutta method; basis augmentation.

\smallskip
\noindent\textbf{MSC 2020.} 65L05, 65L06, 65L20, 65F55.

\section{Introduction}
\label{s4:sec:introduction}
Dynamical low-rank approximation \cite{koch2007} approximates the solution
$A(t)\in\R^{m\times n}$ of a large matrix differential equation
\begin{equation}\label{s4:eq:ode}
 \dot A(t)=F(t,A(t)),\qquad t_0\le t\le T,
\end{equation}
by a matrix of rank $r\ll\min(m,n)$, represented in the factorized form
\begin{equation}\label{s4:eq:intro-factorization}
 Y(t)=U(t)S(t)V(t)^\top.
\end{equation}
Here $U\in\R^{m\times r}$ and $V\in\R^{n\times r}$ have orthonormal
columns, and $S\in\R^{r\times r}$ is nonsingular. Evolving these factors
allows the column and row spaces to follow the solution, while avoiding
the storage of a full $m\times n$ matrix. This approach has been used for
many kinetic equations; see, e.g.,
\cite{vlasov2018,einkemmer2023alfven,boltzmann2025,vpfp2026,kineticreview2025}
and the references therein.

The rank constraint is imposed by projecting the right-hand side of
\eqref{s4:eq:ode} onto the tangent space of the rank-$r$ manifold:
\begin{equation}\label{s4:eq:intro-projected}
 \dot Y(t)=P_r(Y(t))F(t,Y(t)).
\end{equation}
This selects, among the velocities tangent to the manifold, the one
closest to $F(t,Y)$ in the Frobenius norm. At $Y=USV^\top$, the
orthogonal projector is
\begin{equation}\label{s4:eq:intro-projector}
 P_r(Y)Z=ZVV^\top-UU^\top ZVV^\top+UU^\top Z,
 \qquad Z\in\R^{m\times n}.
\end{equation}
Equation~\eqref{s4:eq:intro-projected} can be written as differential
equations for $U$, $S$, and $V$, but the equations for $U$ and $V$
involve $S^{-1}$ \cite{koch2007}. Since the nonzero singular values of
$Y$ are those of $S$, direct time integration of these equations may
require very small step sizes when the approximation has small singular
values. Such values need not indicate a poor approximation: they also
occur when the chosen rank exceeds the rank needed to represent the
solution accurately. Thus a useful low-rank integrator should remain
accurate without step-size restrictions caused by small singular values.

The projector-splitting integrator of Lubich and Oseledets
\cite{lubich2014} avoids the inverse of $S$ by splitting
\eqref{s4:eq:intro-projected} according to the three terms of
\eqref{s4:eq:intro-projector}. Kieri, Lubich, and Walach
\cite{kieri2016} proved a first-order error bound whose constants and
step-size restriction are independent of small singular values. This
property is referred to as \emph{robustness}. The negative middle term
in \eqref{s4:eq:intro-projector}, however, gives a backward substep
that can cause instability for dissipative problems \cite{stability2026}.
The basis-update and Galerkin (BUG) integrator of Ceruti and Lubich
\cite{bug2022} avoids this backward substep. It first updates the column
and row bases independently and then evolves the small core by a
Galerkin equation in the new bases. The method also has a robust
first-order error bound.

Obtaining higher order while retaining robustness requires further
analysis. The Strang composition of projector splitting has classical
order two, but its second-order error bound is not generally uniform in
small singular values \cite{strang2026}. Robust second-order BUG
integrators have been developed using the midpoint rule
\cite{midpoint2024} and parallel basis updates \cite{parallel2025}.
Nobile and Riffaud \cite{nobile2026} develop BUG integrators for general
explicit Runge--Kutta schemes, with error bounds that retain the
underlying time-discretization order up to projection and truncation
errors. Hauck, Kusch, and Schotth\"ofer \cite{hauck2026}
establish general-order BUG methods that also allow implicit or
exponential integration of the Galerkin equation. These high-order BUG
methods use basis augmentation followed by rank truncation. Other
approaches include projected Runge--Kutta methods \cite{kieri2019},
which combine tangent-space projection with low-rank approximation,
and randomized low-rank Runge--Kutta methods \cite{randomized2025}.

In this paper, we develop and analyze robust high-order BUG integrators
without basis augmentation or rank truncation. We consider the BUG
update with a prescribed increment $\Delta A$ approximating the change
in the solution over one time step. The same increment is used in the
basis and core updates, and both bases retain $r$ columns. Our main
result extends the relaxed-exactness analysis of the projector-splitting
integrator in \cite{matrix} to the two-sided projection in BUG. It
provides an accuracy guarantee for general increments, from which
high-order time discretizations can be obtained.

Specifically, let $Y_0$ and $Y_1$ be the input and output of a BUG step.
For a target matrix $A_*$ satisfying the column and row range conditions
stated in Section~\ref{s4:sec:static}, Theorem~\ref{s4:thm:static} gives
\begin{equation}\label{s4:eq:intro-relaxed}
 \|Y_1-A_*\|\le C\|Y_0+\Delta A-A_*\|,
\end{equation}
where $\|\cdot\|$ denotes the Frobenius norm. The constant $C$ depends
on the alignment of the target spaces with the old bases, rather than
on the magnitudes of the nonzero singular values of $A_*$. This
\emph{relaxed exactness} extends the exactness property of BUG
\cite{bug2022} to approximate increments: the error in the low-rank
output is controlled by the error in the prescribed increment.

For the matrix differential equation \eqref{s4:eq:ode}, we obtain robust
local and global error bounds under the boundedness, Lipschitz, and
normal-component assumptions of \cite{matrix}. In particular, if each
increment approximates the exact solution change from the current
numerical state with uniform error $O(h^{p+1}+h\varepsilon_r)$, then
Theorem~\ref{s4:thm:global} yields, for $p\ge1$,
\begin{equation}\label{s4:eq:intro-global}
 \max_{t_0\le t_k\le T}\|Y_k-A(t_k)\|
 \le C(\delta+\varepsilon_r+h^p),\qquad 0<h\le h_0.
\end{equation}
Here $h$ is the time step, $\delta$ bounds the initial error, and
$\varepsilon_r$ bounds the normal component of $F$ at matrices of rank
at most $r$, as specified in \eqref{s4:eq:normal}. The constants $C$ and
$h_0$ are independent of $\delta$, $\varepsilon_r$, and small singular
values. The estimates also cover approximations whose actual rank is
less than $r$.

As concrete high-order methods, we use explicit Runge--Kutta increments
and compute every internal stage and final approximation by the BUG
update. Under the sufficient coefficient conditions established in
\cite{matrix}, the resulting methods have robust orders $p=1,2,3,4$
for forward Euler, explicit midpoint, Heun's third-order method, and
Merson's fourth-order method, respectively
(Corollary~\ref{s4:cor:examples}). Every stage and final approximation
has rank at most $r$, and no augmentation or truncation is required.
The additional stage conditions distinguish this result from
\cite{nobile2026}, which covers arbitrary explicit Runge--Kutta tableaux
of the prescribed classical order using augmented spaces.

We also study basis augmentation and rank truncation within the same
accuracy framework. Augmentation has useful roles beyond high-order
accuracy. In the rank-adaptive BUG integrator \cite{rankadaptive2022},
including the old bases lets the Galerkin step start from the original
numerical solution without a projection error and supports norm and
energy preservation up to truncation error. Basis directions can also
be selected to retain prescribed moments \cite{conservative2023}.
For the midpoint BUG integrator, truncation back to rank $r$ preserves
the robust second-order local error bound without a separate
truncation-error term \cite{midpoint2024}.

Here, we extend the two-sided relaxed-exactness estimate to augmented
BUG updates followed by best rank-$r$ SVD truncation. The additional
directions are arbitrary, provided the BUG sampling spaces are retained.
For general increments, the
output error and the truncation residual are controlled by the increment
error and the normal-component bound. Under our standing assumptions,
applying this estimate at each internal stage and at the endpoint gives
robust high-order Runge--Kutta convergence without a separate
truncation-error assumption (Corollary~\ref{s4:cor:svd-global}). The
global bound in \cite{nobile2026} retains a separate truncation-defect
parameter. Moreover, for the same input and increment, augmentation
followed by this truncation cannot increase
the approximation error to $Y_0+\Delta A$ when the augmented spaces
contain the unaugmented ones (Proposition~\ref{s4:prop:nested-svd}).
Thus basis directions can be added for further requirements while
retaining the high-order accuracy guarantee.

\noindent\textbf{Paper organization.}
Section~\ref{s4:sec:method} defines the BUG update and proves relaxed
exactness. Section~\ref{s4:sec:error} gives the local and global error
bounds and their Runge--Kutta applications.
Section~\ref{s4:sec:augmentation} treats basis augmentation and rank
truncation. The experiments in Section~\ref{s4:sec:numerics} demonstrate
robust high-order convergence without augmentation and compare BUG with
augmented BUG for kinetic shear Alfv\'en waves.

\section{The BUG integrator and relaxed exactness}
\label{s4:sec:method}
We first present the BUG integrator with a fixed increment and then extend
the relaxed-exactness estimate in \cite{matrix} to this integrator.

\subsection{The BUG integrator with a fixed increment}
\label{s4:sec:definition}
Write $\|\cdot\|$ for the Frobenius norm, $\|\cdot\|_2$ for the spectral norm,
and $\langle\cdot,\cdot\rangle_F$ for the Frobenius inner product.
Fix $1\le r\le\min(m,n)$.

Let
\[
 Y_0=U_0S_0V_0^\top,\qquad U_0^\top U_0=V_0^\top V_0=I_r,
\]
where $U_0\in\R^{m\times r}$, $V_0\in\R^{n\times r}$ and
$S_0\in\R^{r\times r}$. The BUG integrator \cite{bug2022} with a fixed
increment $\Delta A\in\R^{m\times n}$ is defined by
\begin{equation}\label{s4:eq:bug}
 \begin{aligned}
 K&=(Y_0+\Delta A)V_0,& L&=(Y_0+\Delta A)^\top U_0,\\
 S_1&=U_1^\top(Y_0+\Delta A)V_1,& Y_1&=U_1S_1V_1^\top.
 \end{aligned}
\end{equation}
Here $U_1$ and $V_1$ are the orthonormal factors in thin QR factorizations
of $K$ and $L$, respectively, each with $r$ columns.
The two basis updates use the old bases and the same increment and can be
computed in parallel. They give
\begin{equation}\label{s4:eq:sampling}
 \begin{aligned}
 &(I-U_1U_1^\top)(Y_0+\Delta A)V_0=0,\\
 &U_0^\top(Y_0+\Delta A)(I-V_1V_1^\top)=0,\\
 &Y_1=U_1U_1^\top(Y_0+\Delta A)V_1V_1^\top.
 \end{aligned}
\end{equation}

\begin{remark}\label{s4:rem:bug-relation}
The update \eqref{s4:eq:bug} is obtained from the BUG construction
in \cite{bug2022} by using the same increment $\Delta A$ in all three
updates. For the nonlinear equation
\eqref{s4:eq:ode}, independently integrating the $K$, $L$, and $S$
subproblems generally gives different increments and hence a different method.
\end{remark}

\subsection{Relaxed exactness}
\label{s4:sec:static}
To relate the output error to the increment error, let
$A_*\in\R^{m\times n}$ satisfy the range condition in
\cite{matrix} and its transpose,
\[
 \range(A_*)=\range(A_*V_0),\qquad
 \range(A_*^\top)=\range(A_*^\top U_0).
\]
Equivalently, there exist $Q_R\in\R^{r\times n}$ and
$Q_L\in\R^{m\times r}$ such that
\begin{equation}\label{s4:eq:reconstruction}
 A_*=(A_*V_0)Q_R=Q_L(U_0^\top A_*).
\end{equation}

\begin{theorem}[Relaxed exactness]\label{s4:thm:static}
Let $U_0,V_0$ have $r$ orthonormal columns, and let $U_1,V_1$ be any matrices
with orthonormal columns satisfying the two sampling identities in
\eqref{s4:eq:sampling}.
For any $A_*,Q_L,Q_R$ satisfying \eqref{s4:eq:reconstruction}, the output
$Y_1=U_1U_1^\top(Y_0+\Delta A)V_1V_1^\top$ satisfies
\begin{equation}\label{s4:eq:static-bound}
 \|Y_1-A_*\|
 \le\sqrt{\max\{1,\|Q_L\|_2^2+\|Q_R\|_2^2\}}\,
       \|Y_0+\Delta A-A_*\|.
\end{equation}
\end{theorem}
\begin{proof}
Set $P=U_1U_1^\top$, $R=V_1V_1^\top$, and
$E=Y_0+\Delta A-A_*$. By \eqref{s4:eq:sampling} and
\eqref{s4:eq:reconstruction},
\[
 \begin{aligned}
 (I-P)A_*&=(I-P)(Y_0+\Delta A-E)V_0Q_R
             =-(I-P)EV_0Q_R,\\
 PA_*(I-R)&=PQ_LU_0^\top(Y_0+\Delta A-E)(I-R)
             =-PQ_LU_0^\top E(I-R).
 \end{aligned}
\]
Hence
\begin{equation}\label{s4:eq:three-terms}
 \begin{aligned}
 Y_1-A_*&=PER-(I-P)A_*-PA_*(I-R)\\
        &=PER+(I-P)EV_0Q_R+PQ_LU_0^\top E(I-R).
 \end{aligned}
\end{equation}
Since $P(I-P)=R(I-R)=0$, the three terms in the last line are orthogonal
in the Frobenius inner product. Therefore,
\begin{equation}\label{s4:eq:static-squares}
 \|Y_1-A_*\|^2
 \le \|PER\|^2+\|Q_R\|_2^2\|(I-P)E\|^2
                +\|Q_L\|_2^2\|E(I-R)\|^2.
\end{equation}
Splitting the last two norms gives
\[
 \begin{aligned}
 \|Y_1-A_*\|^2
 &\le \|PER\|^2+\|Q_L\|_2^2\|PE(I-R)\|^2
                  +\|Q_R\|_2^2\|(I-P)ER\|^2\\
 &\quad +(\|Q_L\|_2^2+\|Q_R\|_2^2)\|(I-P)E(I-R)\|^2\\
 &\le \max\{1,\|Q_L\|_2^2+\|Q_R\|_2^2\}\,\|E\|^2.
 \end{aligned}
\]
Taking square roots proves \eqref{s4:eq:static-bound}.
\end{proof}

In particular, if $\Delta A=A_*-Y_0$, then $Y_1=A_*$, recovering
the exactness property of the BUG integrator \cite{bug2022} under the stated
range conditions.

\section{Robust error analysis}
\label{s4:sec:error}
We derive robust local and global error bounds and then construct robust
high-order methods based on explicit Runge--Kutta schemes without basis augmentation.

\subsection{A robust local error bound}
\label{s4:sec:local}
\paragraph{Standing assumptions.}
Let $0\le t_0<T$ and
$\mathcal M_{\le r}=\{X\in\R^{m\times n}:\rank(X)\le r\}$.
We adopt the assumptions of \cite{matrix}.
The field $F:[0,T]\times\R^{m\times n}\to\R^{m\times n}$ is continuous.
For constants $B>0$, $L_F>0$, and $\varepsilon_r\ge0$, it satisfies
\begin{gather}
 \|F(t,X)\|\le B,\qquad
 \|F(t,X)-F(t,Z)\|\le L_F\|X-Z\|,
 \qquad X,Z\in\R^{m\times n},\label{s4:eq:BL}\\
 \|(I-XX^+)F(t,X)(I-X^+X)\|_2\le\varepsilon_r,
 \qquad X\in\mathcal M_{\le r},\label{s4:eq:normal}
\end{gather}
for all $t\in[0,T]$. Here $X^+$ is the Moore--Penrose pseudoinverse.
The matrices $XX^+$ and $X^+X$ are the orthogonal projectors onto
$\range(X)$ and $\range(X^\top)$, respectively.

Define
\begin{equation}\label{s4:eq:functions}
 \varphi(t)=\frac{e^{L_Ft}-1}{L_F},\qquad
 \gamma(t)=\frac1{\sqrt{1-rL_F^2t^2}},\qquad
 h_* =\frac1{\sqrt{2r}\,L_F}.
\end{equation}
Here $\varphi$ is defined for $t\ge0$ and $\gamma$ for
$0\le t<1/(\sqrt r L_F)$. Set
\begin{equation}\label{s4:eq:kappa}
 \gamma_*:=\gamma(h_*)=\sqrt2,\qquad
 \kappa(h):=\sqrt2\gamma(h),\qquad \kappa_*:=\kappa(h_*)=2.
\end{equation}

To apply relaxed exactness to the differential equation, we construct a
nearby matrix satisfying both range conditions.

\begin{lemma}\label{s4:lem:comparison}
Under \eqref{s4:eq:BL}--\eqref{s4:eq:normal}, let $A$ solve
\eqref{s4:eq:ode} with $A(t_0)=Y_0=U_0S_0V_0^\top$ as in
\eqref{s4:eq:bug}. For $0<h\le h_*$ and $t_1=t_0+h\le T$, there exist
$\widetilde A_1\in\mathcal M_{\le r}$, $Q_L\in\R^{m\times r}$, and
$Q_R\in\R^{r\times n}$ such that
\begin{equation}\label{s4:eq:comparison-error}
 \|A(t_1)-\widetilde A_1\|\le\sqrt{2\min(m,n)}\,\varepsilon_r\varphi(h),
\end{equation}
and the range conditions hold in the form
\begin{equation}\label{s4:eq:comparison-reconstruction}
 \widetilde A_1=Q_L(U_0^\top\widetilde A_1)
               =(\widetilde A_1V_0)Q_R,
\end{equation}
with
\begin{equation}\label{s4:eq:comparison-factors}
 U_0^\top Q_L=Q_RV_0=I_r,\qquad
 \max\{\|Q_L\|_2,\|Q_R\|_2\}\le\gamma(h).
\end{equation}
\end{lemma}
\begin{proof}
The comparison-matrix lemma in \cite[Lemma~3.2]{matrix} gives $X_R,Q_R$ such that
\begin{equation}\label{s4:eq:right-comparison}
 \begin{gathered}
 \|A(t_1)-X_R\|\le\sqrt{\min(m,n)}\,\varepsilon_r\varphi(h),\\
 X_R=(X_RV_0)Q_R,\qquad Q_RV_0=I_r,\qquad \|Q_R\|_2\le\gamma(h).
 \end{gathered}
\end{equation}
The transposed equation $\dot Z=F(t,Z^\top)^\top$, $Z(t_0)=Y_0^\top$,
satisfies the same assumptions with unchanged constants. Applying the lemma
to this equation gives
\begin{equation}\label{s4:eq:left-comparison}
 \begin{gathered}
 \|A(t_1)-X_L\|\le\sqrt{\min(m,n)}\,\varepsilon_r\varphi(h),\\
 X_L=Q_L(U_0^\top X_L),\qquad U_0^\top Q_L=I_r,
 \qquad \|Q_L\|_2\le\gamma(h).
 \end{gathered}
\end{equation}
To combine these two estimates, let $P_c,R_c$ project onto
$\range(Q_L)$ and $\range(Q_R^\top)$, respectively,
and set
\begin{equation}\label{s4:eq:two-sided-comparison}
 \widetilde A_1=P_cA(t_1)R_c.
\end{equation}
Then $\widetilde A_1\in\mathcal M_{\le r}$. Since $P_cX_L=X_L$ and
$X_RR_c=X_R$,
\[
 \|(I-P_c)A(t_1)\|\le\sqrt{\min(m,n)}\,\varepsilon_r\varphi(h),\qquad
 \|A(t_1)(I-R_c)\|\le\sqrt{\min(m,n)}\,\varepsilon_r\varphi(h).
\]
The two terms in
\[
 A(t_1)-\widetilde A_1=(I-P_c)A(t_1)+P_cA(t_1)(I-R_c)
\]
are orthogonal, which gives \eqref{s4:eq:comparison-error}.
The normalizations in \eqref{s4:eq:right-comparison} and
\eqref{s4:eq:left-comparison} imply
\[
 Q_LU_0^\top P_c=P_c,\qquad R_cV_0Q_R=R_c.
\]
Thus $\widetilde A_1$ satisfies the range conditions
\eqref{s4:eq:comparison-reconstruction}.
\end{proof}

Combining this comparison with Theorem~\ref{s4:thm:static} gives the local
error bound.

\begin{theorem}\label{s4:thm:local}
Under \eqref{s4:eq:BL}--\eqref{s4:eq:normal}, let $A$ solve
\eqref{s4:eq:ode} with $A(t_0)=Y_0=U_0S_0V_0^\top$ as in
\eqref{s4:eq:bug}. If $0<h\le h_*$ and $t_1=t_0+h\le T$, then for every
increment $\Delta A$, the BUG output satisfies
\begin{equation}\label{s4:eq:local}
 \|Y_1-A(t_1)\|
 \le \kappa(h)\|Y_0+\Delta A-A(t_1)\|
       +(1+\kappa(h))\sqrt{2\min(m,n)}\,\varepsilon_r\varphi(h).
\end{equation}
\end{theorem}
\begin{proof}
Let $\widetilde A_1,Q_L,Q_R$ be given by Lemma~\ref{s4:lem:comparison}.
Since $\|Q_L\|_2,\|Q_R\|_2\le\gamma(h)$ and $\gamma(h)\ge1$, the coefficient
in Theorem~\ref{s4:thm:static} is at most $\sqrt2\gamma(h)=\kappa(h)$.
Applying the theorem with $A_*=\widetilde A_1$ gives
\[
 \begin{aligned}
 \|Y_1-A(t_1)\|
 &\le\kappa(h)\|Y_0+\Delta A-\widetilde A_1\|+\|\widetilde A_1-A(t_1)\|\\
 &\le\kappa(h)\|Y_0+\Delta A-A(t_1)\|
       +(1+\kappa(h))\|\widetilde A_1-A(t_1)\|.
 \end{aligned}
\]
Together with \eqref{s4:eq:comparison-error}, this proves \eqref{s4:eq:local}.
\end{proof}

\subsection{Global error for general increments}
\label{s4:sec:global}
Let $t_k=t_0+kh$, and let $Y_k$ be the numerical approximation with
$\rank(Y_k)\le r$. To apply the local bound at each step, define the local
exact solution and increment defect by
\begin{equation}\label{s4:eq:restart}
 \dot A_k(t)=F(t,A_k(t)),\qquad A_k(t_k)=Y_k,
 \qquad t_k\le t\le t_{k+1},
\end{equation}
\begin{equation}\label{s4:eq:defect}
 d_k=\|Y_k+\Delta A_k-A_k(t_{k+1})\|,\qquad
 d_h=\max_{t_{k+1}\le T}d_k.
\end{equation}
The increment $\Delta A_k$ approximates $A_k(t_{k+1})-Y_k$.
In order estimates, constants are uniform in the step number, $h$, the initial
error, $\varepsilon_r$, small singular values, and the choice of the
$r$-column QR factors.

\begin{theorem}\label{s4:thm:global}
Assume \eqref{s4:eq:BL}--\eqref{s4:eq:normal} and $0<h\le h_*$, and apply
\eqref{s4:eq:bug} at each step. Let $A$ solve \eqref{s4:eq:ode} with
$\|Y_0-A(t_0)\|\le\delta$. Then
\begin{equation}\label{s4:eq:global}
 \max_{t_k\le T}\|Y_k-A(t_k)\|
 \le e^{L_FT}\delta+\varphi(T)
       \left(2\frac{d_h}{h}+3\sqrt{2\min(m,n)}\,\varepsilon_r\right).
\end{equation}
Consequently, for an integer $p\ge1$,
\begin{equation}\label{s4:eq:order-transfer}
 d_h=O(h^{p+1}+h\varepsilon_r)
 \quad\Longrightarrow\quad
 \max_{t_k\le T}\|Y_k-A(t_k)\|=O(\delta+\varepsilon_r+h^p).
\end{equation}
\end{theorem}
\begin{proof}
Following \cite[Theorem~4.1]{matrix}, set $e_k=\|Y_k-A(t_k)\|$.
For $t_k\le t\le t_{k+1}$, the integral equations for $A_k$ and $A$ give
\[
 A_k(t)-A(t)=Y_k-A(t_k)
 +\int_{t_k}^{t}\bigl(F(t',A_k(t'))-F(t',A(t'))\bigr)\,dt'.
\]
The Lipschitz bound \eqref{s4:eq:BL} implies
\[
 \|A_k(t)-A(t)\|
 \le e_k+L_F\int_{t_k}^{t}\|A_k(t')-A(t')\|\,dt'.
\]
By Gronwall's inequality,
\[
 \|A_k(t)-A(t)\|\le e^{L_F(t-t_k)}e_k.
\]
Apply Theorem~\ref{s4:thm:local} to $A_k$ on $[t_k,t_{k+1}]$.
Since $\kappa(h)\le\kappa_*=2$, the preceding estimate yields
\begin{equation}\label{s4:eq:global-recursion}
 \begin{aligned}
 e_{k+1}
 &\le\|Y_{k+1}-A_k(t_{k+1})\|
       +\|A_k(t_{k+1})-A(t_{k+1})\|\\
 &\le e^{L_Fh}e_k+\kappa(h)d_k
             +(1+\kappa(h))\sqrt{2\min(m,n)}\,\varepsilon_r\varphi(h)\\
 &\le e^{L_Fh}e_k+2d_h+3\sqrt{2\min(m,n)}\,\varepsilon_r\varphi(h).
 \end{aligned}
\end{equation}
Starting from $e_0\le\delta$, iteration gives
\[
 e_k\le e^{L_Fkh}\delta+
 \bigl(2d_h+3\sqrt{2\min(m,n)}\,\varepsilon_r\varphi(h)\bigr)
       \sum_{j=0}^{k-1}e^{L_Fjh}.
\]
To estimate the sum, use
\[
 \varphi(h)\ge h,\qquad
 \varphi(h)\sum_{j=0}^{k-1}e^{L_Fjh}
 =\frac{e^{L_Fkh}-1}{L_F}=\varphi(kh).
\]
Since $kh=t_k-t_0\le T$, we obtain
\[
 \begin{aligned}
 e_k
 &\le e^{L_Fkh}\delta+\varphi(kh)
       \left(\frac{2d_h}{\varphi(h)}
             +3\sqrt{2\min(m,n)}\,\varepsilon_r\right)\\
 &\le e^{L_FT}\delta+\varphi(T)
       \left(2\frac{d_h}{h}
             +3\sqrt{2\min(m,n)}\,\varepsilon_r\right).
 \end{aligned}
\]
This proves \eqref{s4:eq:global}. Substituting
$d_h=O(h^{p+1}+h\varepsilon_r)$ gives \eqref{s4:eq:order-transfer}.
\end{proof}

\subsection{Runge--Kutta methods}
\label{s4:sec:rk}
To construct increments with the accuracy required by the global error bound,
we consider an $s$-stage explicit Runge--Kutta method, $s\ge1$, with Butcher tableau
\[
\begin{array}{c|ccccc}
0      & 0      & 0      & \cdots & 0      & 0\\
c_2    & a_{21} & 0      & \cdots & 0      & 0\\
c_3    & a_{31} & a_{32} & \ddots & \vdots & \vdots\\
\vdots & \vdots & \vdots & \ddots & 0      & 0\\
c_s    & a_{s1} & a_{s2} & \cdots & a_{s,s-1} & 0\\
\hline
       & b_1    & b_2    & \cdots & b_{s-1} & b_s
\end{array}
\]
where $a_{ij}=0$ for $j\ge i$ and $c_i=\sum_{j<i}a_{ij}$.
At time $t_k$, retain a factorization
$Y_k=U_kS_kV_k^\top$ of width $r$. Set $Y_{k,1}=Y_k$ and, successively, form
\begin{equation}\label{s4:eq:rk-increments}
 \begin{aligned}
 t_{k,i}&=t_k+c_i h,& F_{k,i}&=F(t_{k,i},Y_{k,i}),\\
 \Delta A_{k,i}&=h\sum_{j<i}a_{ij}F_{k,j} \quad (i\ge2),&
 \Delta A_k&=h\sum_{i=1}^s b_iF_{k,i}.
 \end{aligned}
\end{equation}
For each $i\ge2$, compute $Y_{k,i}$ by \eqref{s4:eq:bug} with input $Y_k$,
old bases $U_k,V_k$, and increment $\Delta A_{k,i}$. After evaluating all
stages, use the same input and old bases with increment $\Delta A_k$ to compute
$Y_{k+1}$. Algorithm~\ref{s4:alg:bug-rk} summarizes one step.

\begin{algorithm}[htbp]
\caption{One step of the BUG Runge--Kutta method.}
\label{s4:alg:bug-rk}
\begin{algorithmic}[1]
\Require $t_k,h$, factors $Y_k=U_kS_kV_k^\top$ of width $r$, and an explicit
$s$-stage tableau $(a_{ij},b_i,c_i)$
\State $Y_{k,1}\gets Y_k$, $F_{k,1}\gets F(t_k,Y_k)$.
\For{$i=2,\ldots,s$}
 \State $\Delta A_{k,i}\gets h\sum_{j=1}^{i-1}a_{ij}F_{k,j}$.
 \State $K\gets(Y_k+\Delta A_{k,i})V_k$, $L\gets(Y_k+\Delta A_{k,i})^\top U_k$.
 \State Compute the $r$-column orthonormal QR factors $U_{k,i},V_{k,i}$ of $K,L$.
 \State $S_{k,i}\gets U_{k,i}^\top(Y_k+\Delta A_{k,i})V_{k,i}$,
 $Y_{k,i}\gets U_{k,i}S_{k,i}V_{k,i}^\top$.
 \State $F_{k,i}\gets F(t_k+c_i h,Y_{k,i})$.
\EndFor
\State $\Delta A_k\gets h\sum_{i=1}^s b_iF_{k,i}$.
\State $K\gets(Y_k+\Delta A_k)V_k$, $L\gets(Y_k+\Delta A_k)^\top U_k$.
\State Compute the $r$-column orthonormal QR factors $U_{k+1},V_{k+1}$ of $K,L$.
\State $S_{k+1}\gets U_{k+1}^\top(Y_k+\Delta A_k)V_{k+1}$.
\Ensure $Y_{k+1}=U_{k+1}S_{k+1}V_{k+1}^\top$.
\end{algorithmic}
\end{algorithm}
\FloatBarrier

To estimate the increment error, assume for an integer $p\ge1$ that
\begin{equation}\label{s4:eq:regularity}
 \begin{gathered}
 g_k(t):=F(t,A_k(t))\in C^p([t_k,t_{k+1}]),\\
 \sup_{t_k\le t\le t_{k+1}}\|g_k^{(a)}(t)\|\le M_a,
 \qquad 1\le a\le p.
 \end{gathered}
\end{equation}
The bounds hold uniformly for the local solutions $A_k$ arising from all
choices of the $r$-column QR factors, with $M_a$ independent of $k$,
$0<h\le h_0:=\min(1,h_*)$, and small singular values.

\begin{corollary}\label{s4:cor:rk}
Under the standing assumptions and \eqref{s4:eq:regularity}, suppose
$0<c_i\le1$ for $i\ge2$. Use the coefficient conditions of \cite[Theorem~5.1]{matrix}:
\begin{align}
 \sum_{i=1}^s b_i&=1,&
 \sum_{i=1}^s b_i c_i^l&=\frac1{l+1},&&1\le l\le p-1.
 \label{s4:eq:output-moments}
\end{align}
Suppose either $p=s=1$, or $p\ge2$ and integers
$q_i\in\{1,\ldots,p-1\}$, $2\le i\le s$, satisfy
\begin{align}
 \sum_{j<i}a_{ij}c_j^l&=\frac{c_i^{l+1}}{l+1},
 &&1\le l<q_i,\quad 2\le i\le s,\label{s4:eq:row-moments}\\
 a_{ij}\ne0&\ \Longrightarrow\ q_j\ge q_i-1,
 &&2\le j<i\le s,\label{s4:eq:dependencies}\\
 b_i\ne0&\ \Longrightarrow\ q_i=p-1,
 &&2\le i\le s.\label{s4:eq:support}
\end{align}
Then, for $0<h\le h_0$, the BUG stages and final increment satisfy
\begin{align}
 \|Y_{k,i}-A_k(t_{k,i})\|&=O(h^{q_i+1}+h\varepsilon_r),
 &&2\le i\le s,\label{s4:eq:stage-error}\\
 d_k&=O(h^{p+1}+h^2\varepsilon_r),\label{s4:eq:rk-defect}\\
 \|Y_{k+1}-A_k(t_{k+1})\|&=O(h^{p+1}+h\varepsilon_r).
 \label{s4:eq:rk-local}
\end{align}
With $\|Y_0-A(t_0)\|\le\delta$, the global error is
$O(\delta+\varepsilon_r+h^p)$.
The constants may depend on the tableau, $r$, $\min(m,n)$, $L_F$, $T$, $p$,
and $M_1,\ldots,M_p$, but not on small singular values or the choice of QR factors.
\end{corollary}
\begin{proof}
We refer to the proof of \cite[Theorem~5.1]{matrix}, replacing its local
error estimate by Theorem~\ref{s4:thm:local}.
Fix $k$ and put $\eta_i=\|Y_{k,i}-A_k(t_{k,i})\|$; then $\eta_1=0$.
For $p\ge2$, the Taylor-quadrature estimate there gives
\begin{equation}\label{s4:eq:stage-increment-bound}
 \|Y_k+\Delta A_{k,i}-A_k(t_{k,i})\|
 \le hL_F\sum_{j<i}|a_{ij}|\eta_j+O(h^{q_i+1}).
\end{equation}
Apply Theorem~\ref{s4:thm:local} with duration $c_i h$ and input
$Y_k,U_k,V_k$.
Since $0<c_i\le1$, we have $\kappa(c_i h)\le2$ and
$\varphi(c_i h)\le h e^{L_F}$ for $h\le1$, which gives
\begin{equation}\label{s4:eq:stage-recursion}
 \eta_i\le 2hL_F\sum_{j<i}|a_{ij}|\eta_j
      +O(h^{q_i+1})+3\sqrt{2\min(m,n)}\,e^{L_F}h\varepsilon_r.
\end{equation}
Only the constants in this recursion differ from those in the cited proof.
Its stage induction therefore gives \eqref{s4:eq:stage-error}.
For every $i\ge2$ with $b_i\ne0$, condition~\eqref{s4:eq:support} ensures
that the stage error contributes
\[
 h\,O(h^{q_i+1}+h\varepsilon_r)=O(h^{p+1}+h^2\varepsilon_r)
\]
to the final increment. The final quadrature argument in the cited proof
then gives \eqref{s4:eq:rk-defect}.
For $p=s=1$, the same quadrature argument gives $d_k=O(h^2)$.
The local and global conclusions follow from
Theorems~\ref{s4:thm:local} and \ref{s4:thm:global}.
\end{proof}

\begin{corollary}\label{s4:cor:examples}
With the corresponding regularity order in \eqref{s4:eq:regularity}, the four
tableaux in Table~\ref{s4:tab:tableaux} give robust orders one, two, three,
and four, respectively. The stage assignments $(q_2,\ldots,q_s)$ are
$(1)$ for RK2-Midpoint, $(1,2)$ for RK3-Heun, and $(1,2,3,3)$ for RK4-Merson.
\end{corollary}
\begin{proof}
The coefficient conditions are verified in \cite[Corollary~5.2]{matrix}.
The result follows from Corollary~\ref{s4:cor:rk}.
\end{proof}

\begin{table}[htbp]
\centering
\[
\begin{array}{c|c}
0&0\\ \hline
&1
\end{array}
\qquad
\begin{array}{c|cc}
0&0&0\\
\frac12&\frac12&0\\ \hline
&0&1
\end{array}
\qquad
\begin{array}{c|ccc}
0&0&0&0\\
\frac13&\frac13&0&0\\
\frac23&0&\frac23&0\\ \hline
&\frac14&0&\frac34
\end{array}
\qquad
\begin{array}{c|ccccc}
0&0&0&0&0&0\\
\frac13&\frac13&0&0&0&0\\
\frac13&\frac16&\frac16&0&0&0\\
\frac12&\frac18&0&\frac38&0&0\\
1&\frac12&0&-\frac32&2&0\\ \hline
&\frac16&0&0&\frac23&\frac16
\end{array}
\]
\caption{Butcher tableaux for Euler, RK2-Midpoint, RK3-Heun, and RK4-Merson,
respectively.}
\label{s4:tab:tableaux}
\end{table}
\FloatBarrier

\paragraph{Implementation and cost.}
\label{s5:sec:cost}
The update \eqref{s4:eq:bug} requires the products $\Delta A V_0$,
$\Delta A^\top U_0$, and $U_1^\top\Delta A V_1$. Contributions from $Y_0$
are computed from its factors; for example, its contribution to the new core
is $(U_1^\top U_0)S_0(V_0^\top V_1)$. At RK stage $i$, the new core
requires $U_{k,i}^\top F_{k,j}V_{k,i}$ for previous stages $j<i$.
If $F_{k,j}$ is stored in factored form as
$F_{k,j}=\sum_\ell X_\ell Z_\ell^\top$, these products are computed by
\[
 U_{k,i}^\top F_{k,j}V_{k,i}
 =\sum_\ell (U_{k,i}^\top X_\ell)(Z_\ell^\top V_{k,i}).
\]
The stage field factors are retained throughout the step and contracted
with the new bases at each stage and at the endpoint.

Once these increment products are available, one reconstruction uses two
thin QR factorizations and operations on the factors, requiring
$O((m+n)r^2+r^3)$ work and $O((m+n)r+r^2)$ working storage. It requires no
singular value decomposition. An $s$-stage step has $s$ field evaluations and
$s$ reconstructions: $s-1$ for the internal stages after the first, and one
for the endpoint. Thus its factorization and remaining factor operations
cost $O(s((m+n)r^2+r^3))$.

The additional work and storage for the field values depend on the problem
and their representation. Direct dense assembly of all RK sums and their
projections costs $O(s^2mn+smnr)$ per step, in addition to the field
evaluations, and retaining the $s$ field matrices uses $O(smn)$ storage.
For factored field representations, the work and storage depend on the
number and widths of the factors and the resulting contractions.
Section~\ref{s4:sec:augmentation} allows
augmentation at the stages and endpoint, followed by small-core SVD
truncation. The temporary widths depend on the added directions, while
the matrices used for field evaluations and the next step have rank at most $r$.
\FloatBarrier

\section{Basis augmentation and rank truncation}
\label{s4:sec:augmentation}
We extend the accuracy analysis to augmented bases and show that ordinary
rank truncation preserves robust high-order convergence.

The preceding results establish high-order accuracy without augmentation.
Additional basis directions can nevertheless be useful for structure
preservation. In the original BUG integrator \cite{bug2022}, the Galerkin
step in the new bases $U_1,V_1$ has initial core $U_1^\top Y_0V_1$, which
represents
\[
 U_1U_1^\top Y_0V_1V_1^\top.
\]
This projection need not equal $Y_0$. Even if the Galerkin evolution
preserves the Frobenius norm, it preserves the norm of this projected
starting value, which may be smaller than $\|Y_0\|$. The augmented BUG
integrator includes the old bases $U_0,V_0$ so that its Galerkin step
starts from $Y_0$ \cite{rankadaptive2022}. Other directions can be chosen
to represent prescribed moments \cite{conservative2023}.

We allow arbitrary additional directions while retaining the two sampling
spaces of the BUG update. For a general increment, we bound both the
error after best rank-$r$ SVD truncation and the truncation residual by
the increment error and the normal-component bound. Applying these
estimates at the internal stages and the endpoint then gives robust
high-order convergence without a separate truncation-error assumption.

\subsection{The augmented update}
\label{s4:sec:augmented-update}
Let $U_{\mathrm{add}}$ and $V_{\mathrm{add}}$ contain the additional column
and row directions; either may have no columns. Algorithm~\ref{s4:alg:augmentation}
augments the sampling spaces, projects the increment target, and truncates
the small core. Here $\operatorname{orth}$ denotes an orthonormal basis of
the column space, and $\mathcal T_r$ denotes a best approximation of rank
at most $r$ in the Frobenius norm, obtained by SVD.

\begin{algorithm}[htbp]
\caption{BUG update with basis augmentation and rank truncation.}
\label{s4:alg:augmentation}
\begin{algorithmic}[1]
\Require Factors $Y_0=U_0S_0V_0^\top$ of width $r$, increment $\Delta A$, and
additional directions $U_{\mathrm{add}},V_{\mathrm{add}}$
\State $K\gets(Y_0+\Delta A)V_0$, $L\gets(Y_0+\Delta A)^\top U_0$.
\State $\widehat U_1\gets\operatorname{orth}[K,U_{\mathrm{add}}]$,
$\widehat V_1\gets\operatorname{orth}[L,V_{\mathrm{add}}]$.
\State $\widehat S_1\gets\widehat U_1^\top(Y_0+\Delta A)\widehat V_1$.
\State Define $\widehat Y_1=\widehat U_1\widehat S_1\widehat V_1^\top$.
\State Compute a best rank-$r$ SVD approximation
$\Phi\Sigma\Psi^\top=\mathcal T_r(\widehat S_1)$.
\State $U_1\gets\widehat U_1\Phi$, $S_1\gets\Sigma$,
$V_1\gets\widehat V_1\Psi$.
\Ensure $\widehat Y_1$ and $Y_1=U_1S_1V_1^\top=\mathcal T_r(\widehat Y_1)$,
represented by their factors.
\end{algorithmic}
\end{algorithm}
\FloatBarrier

Since the outer factors have orthonormal columns, this computes
\begin{equation}\label{s4:eq:augmented-truncated}
 \widehat Y_1=\widehat U_1\widehat U_1^\top(Y_0+\Delta A)
                  \widehat V_1\widehat V_1^\top,\qquad
 Y_1=\mathcal T_r(\widehat Y_1).
\end{equation}
Before the next step, represent $Y_1$ by factors of width $r$ as in
\eqref{s4:eq:bug}.

Taking $U_{\mathrm{add}}=U_0$ and $V_{\mathrm{add}}=V_0$ retains the
representation of $Y_0$ in the augmented spaces. To retain prescribed
moments before truncation, a matrix $W\in\R^{n\times d_W}$ specifies the
choice
\begin{equation}\label{s4:eq:moment-spaces}
 U_{\mathrm{add}}=(Y_0+\Delta A)W,\qquad V_{\mathrm{add}}=W.
\end{equation}
It gives
\begin{equation}\label{s4:eq:moment-exactness}
 \widehat Y_1W
 =\widehat U_1\widehat U_1^\top(Y_0+\Delta A)\widehat V_1\widehat V_1^\top W
 =(Y_0+\Delta A)W.
\end{equation}
After truncation, the moment defect satisfies
\[
 \|Y_1W-(Y_0+\Delta A)W\|
 =\|(Y_1-\widehat Y_1)W\|
 \le\|Y_1-\widehat Y_1\|\,\|W\|_2.
\]

\subsection{Accuracy after rank truncation}
\label{s4:sec:svd-accuracy}
The proof of Theorem~\ref{s4:thm:static} uses the new bases through their
orthogonality and the sampling identities. Both properties hold for the
augmented bases, so we first obtain the same bound before truncation.

\begin{corollary}\label{s4:cor:augmentation}
Let $A_*,Q_L,Q_R$ satisfy \eqref{s4:eq:reconstruction}. For the augmented
endpoint in \eqref{s4:eq:augmented-truncated},
\[
 \|\widehat Y_1-A_*\|
 \le\sqrt{\max\{1,\|Q_L\|_2^2+\|Q_R\|_2^2\}}\,
       \|Y_0+\Delta A-A_*\|.
\]
Under the hypotheses of Theorem~\ref{s4:thm:local}, the local bound
\eqref{s4:eq:local} also holds for $\widehat Y_1$.
\end{corollary}
\begin{proof}
The augmented bases satisfy the two sampling identities in
\eqref{s4:eq:sampling}. Apply Theorem~\ref{s4:thm:static}, followed by
the comparison argument in Theorem~\ref{s4:thm:local}.
\end{proof}

\begin{remark}\label{s4:rem:old-space}
If $\widehat U_1\widehat U_1^\top U_0=U_0$, the coefficient before
truncation improves to $\max\{1,\|Q_L\|_2,\|Q_R\|_2\}$, and the local
bound holds with $\kappa(h)$ replaced by $\gamma(h)$.
Indeed, in the notation of the static proof,
$U_0^\top E(I-R)=U_0^\top PE(I-R)$.
The last norm in \eqref{s4:eq:static-squares} can therefore be replaced
by $\|PE(I-R)\|$, leaving three orthogonal error blocks.
\end{remark}

To include the truncation, the projected reference matrix provides a
rank-$r$ competitor for the SVD. This gives a relaxed-exactness estimate
for the complete update.

\begin{proposition}\label{s4:prop:svd-static}
For any $A_*,Q_L,Q_R$ satisfying \eqref{s4:eq:reconstruction}, the output
of Algorithm~\ref{s4:alg:augmentation} satisfies
\begin{equation}\label{s4:eq:svd-static}
 \|Y_1-A_*\|
 \le\sqrt{\max\{4,\|Q_L\|_2^2+\|Q_R\|_2^2\}}\,
       \|Y_0+\Delta A-A_*\|.
\end{equation}
\end{proposition}
\begin{proof}
Set $P=\widehat U_1\widehat U_1^\top$,
$R=\widehat V_1\widehat V_1^\top$, and $E=Y_0+\Delta A-A_*$.
By \eqref{s4:eq:reconstruction},
$\rank(A_*)=\rank((A_*V_0)Q_R)\le r$. Hence
$\rank(PA_*R)\le r$, and SVD optimality gives
\begin{equation}\label{s4:eq:svd-tail-static}
 \|\widehat Y_1-Y_1\|
 \le\|\widehat Y_1-PA_*R\|=\|PER\|.
\end{equation}
Subtracting the truncation residual from \eqref{s4:eq:three-terms} yields
\[
 Y_1-A_*=
 \bigl(PER-(\widehat Y_1-Y_1)\bigr)
 +(I-P)EV_0Q_R+PQ_LU_0^\top E(I-R).
\]
The three terms are orthogonal: the first lies in the product space
defined by $P,R$, the second has its columns in $\range(I-P)$, and the
third has its columns in $\range(P)$ and rows in $\range(I-R)$.
Using \eqref{s4:eq:svd-tail-static} and then splitting $E$ into its four
orthogonal blocks gives
\[
 \begin{aligned}
 \|Y_1-A_*\|^2
 &\le4\|PER\|^2+\|Q_R\|_2^2\|(I-P)E\|^2
                   +\|Q_L\|_2^2\|E(I-R)\|^2\\
 &=4\|PER\|^2+\|Q_L\|_2^2\|PE(I-R)\|^2
                 +\|Q_R\|_2^2\|(I-P)ER\|^2\\
 &\quad+(\|Q_L\|_2^2+\|Q_R\|_2^2)\|(I-P)E(I-R)\|^2\\
 &\le\max\{4,\|Q_L\|_2^2+\|Q_R\|_2^2\}\,\|E\|^2.
 \end{aligned}
\]
Taking square roots proves the result.
\end{proof}

Applying this estimate to the comparison matrix controls both the
endpoint error and the discarded singular values.

\begin{corollary}\label{s4:cor:svd-local}
Under the hypotheses of Theorem~\ref{s4:thm:local}, Algorithm~\ref{s4:alg:augmentation}
satisfies
\begin{align}
 \|Y_1-A(t_1)\|
 &\le2\|Y_0+\Delta A-A(t_1)\|
       +3\sqrt{2\min(m,n)}\,\varepsilon_r\varphi(h),
       \label{s4:eq:svd-local}\\
 \|\widehat Y_1-Y_1\|
 &\le\|Y_0+\Delta A-A(t_1)\|
       +\sqrt{2\min(m,n)}\,\varepsilon_r\varphi(h).
       \label{s4:eq:svd-tail}
\end{align}
\end{corollary}
\begin{proof}
Use $\widetilde A_1,Q_L,Q_R$ from Lemma~\ref{s4:lem:comparison}.
Since $\|Q_L\|_2,\|Q_R\|_2\le\gamma(h)\le\sqrt2$,
Proposition~\ref{s4:prop:svd-static} gives
\[
 \begin{aligned}
 \|Y_1-A(t_1)\|
 &\le2\|Y_0+\Delta A-\widetilde A_1\|
                  +\|\widetilde A_1-A(t_1)\|\\
 &\le2\|Y_0+\Delta A-A(t_1)\|
                  +3\|\widetilde A_1-A(t_1)\|.
 \end{aligned}
\]
The comparison bound proves \eqref{s4:eq:svd-local}.
For \eqref{s4:eq:svd-tail}, use \eqref{s4:eq:svd-tail-static} with
$A_*=\widetilde A_1$ and the contractivity of the orthogonal projections:
\[
 \|\widehat Y_1-Y_1\|
 \le\|Y_0+\Delta A-\widetilde A_1\|
 \le\|Y_0+\Delta A-A(t_1)\|+\|\widetilde A_1-A(t_1)\|.
\]
\end{proof}

To construct an augmented BUG Runge--Kutta method, we now apply
Algorithm~\ref{s4:alg:augmentation} at the internal stages as well as the
endpoint. Each stage uses the same input $Y_k$ and old bases $U_k,V_k$;
its augmented approximation is truncated to rank at most $r$ before
evaluating $F$. The additional directions may vary between stages and
are chosen from quantities already computed. For example, choosing
$U_{\mathrm{add}}=U_k$ and $V_{\mathrm{add}}=V_k$ at every application
retains the old spaces in all augmented bases.

\begin{algorithm}[htbp]
\caption{One step of the augmented BUG Runge--Kutta method.}
\label{s4:alg:augmented-bug-rk}
\begin{algorithmic}[1]
\Require $t_k,h$, factors $Y_k=U_kS_kV_k^\top$ of width $r$, an explicit
$s$-stage tableau $(a_{ij},b_i,c_i)$, and a rule for choosing additional directions
\State $Y_{k,1}\gets Y_k$, $F_{k,1}\gets F(t_k,Y_k)$.
\For{$i=2,\ldots,s$}
 \State $\Delta A_{k,i}\gets h\sum_{j=1}^{i-1}a_{ij}F_{k,j}$.
 \State Choose $U_{\mathrm{add}},V_{\mathrm{add}}$ for this stage.
 \State Apply Algorithm~\ref{s4:alg:augmentation} with
 $(U_k,S_k,V_k)$, $\Delta A_{k,i}$, and these directions to obtain
 $\widehat Y_{k,i}=\widehat U_{k,i}\widehat S_{k,i}\widehat V_{k,i}^\top$
 and $Y_{k,i}=\mathcal T_r(\widehat Y_{k,i})$.
 \State $F_{k,i}\gets F(t_k+c_i h,Y_{k,i})$.
\EndFor
\State $\Delta A_k\gets h\sum_{i=1}^s b_iF_{k,i}$.
\State Choose $U_{\mathrm{add}},V_{\mathrm{add}}$ for the endpoint.
\State Apply Algorithm~\ref{s4:alg:augmentation} with
 $(U_k,S_k,V_k)$, $\Delta A_k$, and these directions to obtain
 $\widehat Y_{k+1}=\widehat U_{k+1}\widehat S_{k+1}\widehat V_{k+1}^\top$
 and $Y_{k+1}=\mathcal T_r(\widehat Y_{k+1})$; represent $Y_{k+1}$ by factors of width $r$.
\Ensure $Y_{k+1}=U_{k+1}S_{k+1}V_{k+1}^\top$.
\end{algorithmic}
\end{algorithm}
\FloatBarrier

The local estimate \eqref{s4:eq:svd-local} applies to each stage with
duration $c_i h$. This supplies the same stage recursion as before and
also controls the truncation errors.

\begin{corollary}\label{s4:cor:svd-global}
Under the assumptions of Theorem~\ref{s4:thm:global}, replace each BUG
update by Algorithm~\ref{s4:alg:augmentation}. Then the global bound
\eqref{s4:eq:global} and the order-transfer conclusion
\eqref{s4:eq:order-transfer} remain valid.
In particular, consider Algorithm~\ref{s4:alg:augmented-bug-rk} under
the tableau conditions of Corollary~\ref{s4:cor:rk} and the uniform
regularity \eqref{s4:eq:regularity} for the resulting restart flows.
For $0<h\le h_0$, its stage, increment, and local errors satisfy
\eqref{s4:eq:stage-error}--\eqref{s4:eq:rk-local}, and
\begin{equation}\label{s4:eq:svd-rk}
 \max_{t_k\le T}\|Y_k-A(t_k)\|=O(\delta+\varepsilon_r+h^p),
 \qquad
 \|\widehat Y_{k+1}-Y_{k+1}\|=O(h^{p+1}+h\varepsilon_r).
\end{equation}
The internal stage truncations satisfy
\begin{equation}\label{s4:eq:svd-stage-tails}
 \|\widehat Y_{k,i}-Y_{k,i}\|=O(h^{q_i+1}+h\varepsilon_r),
 \qquad 2\le i\le s.
\end{equation}
The constants have the dependencies stated in Corollary~\ref{s4:cor:rk}
and are independent of the added directions, subject to those uniform
regularity bounds.
\end{corollary}
\begin{proof}
Use \eqref{s4:eq:svd-local} in place of Theorem~\ref{s4:thm:local} in the
proofs of Theorem~\ref{s4:thm:global} and Corollary~\ref{s4:cor:rk}.
It gives the last bound in \eqref{s4:eq:global-recursion} and the same
stage recursion \eqref{s4:eq:stage-recursion}. Thus the summation, stage
induction, and quadrature argument give \eqref{s4:eq:global},
\eqref{s4:eq:order-transfer}, and
\eqref{s4:eq:stage-error}--\eqref{s4:eq:rk-local} for the augmented method.

For the additional truncation estimates, \eqref{s4:eq:svd-tail} gives
\[
 \begin{aligned}
 \|\widehat Y_{k,i}-Y_{k,i}\|
 &\le\|Y_k+\Delta A_{k,i}-A_k(t_{k,i})\|
       +\sqrt{2\min(m,n)}\,\varepsilon_r\varphi(c_i h)\\
 &=O(h^{q_i+1}+h\varepsilon_r),\\
 \|\widehat Y_{k+1}-Y_{k+1}\|
 &\le d_k+\sqrt{2\min(m,n)}\,\varepsilon_r\varphi(h)\\
 &=O(h^{p+1}+h\varepsilon_r).
 \end{aligned}
\]
The first estimate uses \eqref{s4:eq:stage-increment-bound},
\eqref{s4:eq:stage-error}, and \eqref{s4:eq:dependencies};
the second uses \eqref{s4:eq:rk-defect}.
\end{proof}

The tableaux in Table~\ref{s4:tab:tableaux} therefore give robust orders
one to four in Algorithm~\ref{s4:alg:augmented-bug-rk}.

We finally compare the augmented and unaugmented outputs for the same
increment. The following proposition shows that, when the augmented spaces
contain the spaces selected by the unaugmented update, augmentation followed
by best rank-$r$ SVD truncation does not increase the error in approximating
$Y_0+\Delta A$.

\begin{proposition}\label{s4:prop:nested-svd}
For the same input factors $Y_0=U_0S_0V_0^\top$ and increment $\Delta A$,
let $Y_1^{\mathrm{BUG}}$ be the output of \eqref{s4:eq:bug}, with selected
bases $U_{\mathrm B},V_{\mathrm B}$, and let $\widehat Y_1,Y_1$ be the outputs
of Algorithm~\ref{s4:alg:augmentation}. If
\[
 \range(U_{\mathrm B})\subseteq\range(\widehat U_1),\qquad
 \range(V_{\mathrm B})\subseteq\range(\widehat V_1),
\]
then
\[
 \|Y_0+\Delta A-Y_1\|
 \le\|Y_0+\Delta A-Y_1^{\mathrm{BUG}}\|.
\]
\end{proposition}
\begin{proof}
By definition,
\[
 Y_1^{\mathrm{BUG}}
 =U_{\mathrm B}U_{\mathrm B}^\top(Y_0+\Delta A)
       V_{\mathrm B}V_{\mathrm B}^\top.
\]
The inclusion of both bases implies
\[
 \widehat U_1\widehat U_1^\top Y_1^{\mathrm{BUG}}
       \widehat V_1\widehat V_1^\top=Y_1^{\mathrm{BUG}}.
\]
Thus $Y_1^{\mathrm{BUG}}$ and the truncated output $Y_1$ can both be written
in the augmented bases:
\[
 \begin{aligned}
 Y_1^{\mathrm{BUG}}
 &=\widehat U_1(\widehat U_1^\top Y_1^{\mathrm{BUG}}\widehat V_1)
       \widehat V_1^\top,\\
 Y_1&=\widehat U_1\mathcal T_r(\widehat S_1)\widehat V_1^\top.
 \end{aligned}
\]
By the definition of $\widehat S_1$ and $\widehat Y_1$,
\[
 \widehat U_1^\top(Y_0+\Delta A-\widehat Y_1)\widehat V_1
 =\widehat U_1^\top(Y_0+\Delta A)\widehat V_1-\widehat S_1=0.
\]
Using the Frobenius inner product, we therefore obtain
\[
 \begin{aligned}
 &\langle Y_0+\Delta A-\widehat Y_1,\widehat Y_1-Y_1\rangle_F=\left\langle
   \widehat U_1^\top(Y_0+\Delta A-\widehat Y_1)\widehat V_1,
   \widehat S_1-\mathcal T_r(\widehat S_1)\right\rangle_F=0,
 \end{aligned}
\]
and, similarly,
\[
 \begin{aligned}
 &\langle Y_0+\Delta A-\widehat Y_1,\widehat Y_1-Y_1^{\mathrm{BUG}}\rangle_F=\left\langle
   \widehat U_1^\top(Y_0+\Delta A-\widehat Y_1)\widehat V_1,
   \widehat S_1-\widehat U_1^\top Y_1^{\mathrm{BUG}}\widehat V_1
   \right\rangle_F=0.
 \end{aligned}
\]
Moreover, $Y_1=\mathcal T_r(\widehat Y_1)$ and $\rank(Y_1^{\mathrm{BUG}})\le r$,
so SVD optimality gives
\[
 \|\widehat Y_1-Y_1\|^2
 =\min_{\rank(Z)\le r}\|\widehat Y_1-Z\|^2
 \le\|\widehat Y_1-Y_1^{\mathrm{BUG}}\|^2.
\]
The two zero inner products give the equalities below, while the SVD
bound gives the inequality:
\[
 \begin{aligned}
 \|Y_0+\Delta A-Y_1\|^2
 &=\|Y_0+\Delta A-\widehat Y_1\|^2+\|\widehat Y_1-Y_1\|^2\\
 &\le\|Y_0+\Delta A-\widehat Y_1\|^2
                  +\|\widehat Y_1-Y_1^{\mathrm{BUG}}\|^2\\
 &=\|Y_0+\Delta A-Y_1^{\mathrm{BUG}}\|^2.
 \end{aligned}
\]
\end{proof}

\section{Numerical experiments}
\label{s4:sec:numerics}
We demonstrate robust high-order convergence of the BUG integrator without
augmentation and compare its performance with augmented BUG in two examples.
The augmented methods add the old bases at
every internal stage and at the endpoint, followed by best rank-$r$ SVD
truncation, as in Algorithm~\ref{s4:alg:augmented-bug-rk}.
Reference solutions are computed with DOP853 \cite{hairer1993}, starting from the same
initial matrix as the numerical methods.

\subsection{Robustness to small singular values}
\label{s4:sec:numerical-order}
We first examine the convergence of BUG as the smallest singular value
tends to zero.

\subsubsection{Zero normal component}
\label{s5:sec:zero-normal}
For $m=n=20$ and $r=4$, consider
\begin{equation}\label{s4:eq:numerical-field}
 \begin{aligned}
 F(t,Y)&=\Omega(t,Y)Y+YG(t,Y),\\
 \Omega(t,Y)&=
 \frac{O_0+\sin(2t)O_1+0.7\tanh(\ell^\top Yw)J}{\sqrt{1+\|Y\|^2}},\\
 G(t,Y)&=\frac{G_0+\cos(3t)G_1}{\sqrt{1+\|Y\|^2}}.
 \end{aligned}
\end{equation}
The fixed matrices $O_0,O_1,J,G_0,G_1$ are skew-symmetric and
$\ell,w$ are unit vectors. The exact flow preserves singular values, and
$(I-YY^+)F(t,Y)(I-Y^+Y)=0$ for every $Y$, so $\eps_r=0$.
For fixed orthonormal $U_0,V_0\in\R^{20\times4}$, we take
\begin{equation}\label{s4:eq:numerical-initial}
 Y_0=U_0\operatorname{diag}(1,10^{-1},10^{-2},\sigma)V_0^\top,
 \qquad \sigma\in\{10^{-2},10^{-6},10^{-10},0\},
\end{equation}
and integrate to $T=1$. The case $\sigma=0$ has actual rank three and
working rank four. We measure the relative Frobenius error
$\|Y_N-A_{\mathrm{ref}}(T)\|/\|A_{\mathrm{ref}}(T)\|$.

\begin{figure}[htbp]
\centering
\includegraphics[width=\textwidth]{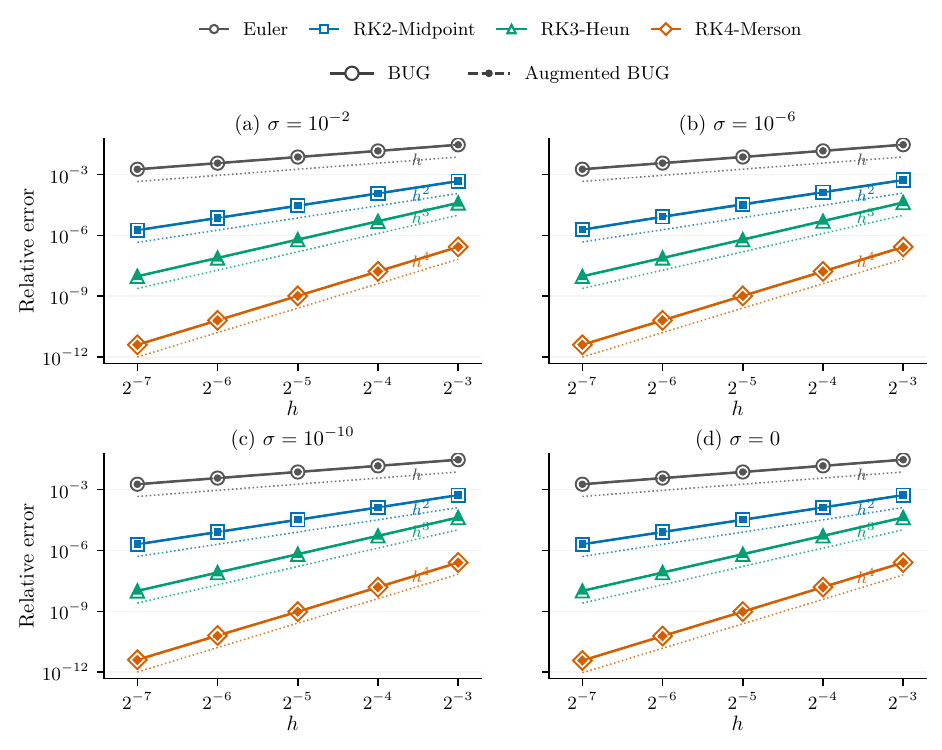}
\caption{Relative errors at $T=1$ for \eqref{s4:eq:numerical-field} with
working rank $r=4$. Each panel corresponds to one value of $\sigma$.
Large open markers and solid lines denote BUG; small filled markers and
dashed lines denote augmented BUG. Dotted lines indicate orders one through four.}
\label{s5:fig:convergence}
\end{figure}

Figure~\ref{s5:fig:convergence} shows first-, second-, third-, and
fourth-order convergence for BUG with Euler, RK2-Midpoint, RK3-Heun, and
RK4-Merson, respectively. These orders persist across all four values of
$\sigma$, including the rank-deficient initial value, and the errors are
nearly identical to those of augmented BUG. Thus the predicted robust
convergence orders are attained without basis augmentation.

\subsubsection{Nonzero normal component}
\label{s5:sec:nonzero-normal}
\begin{figure}[htbp]
\centering
\includegraphics[width=\textwidth]{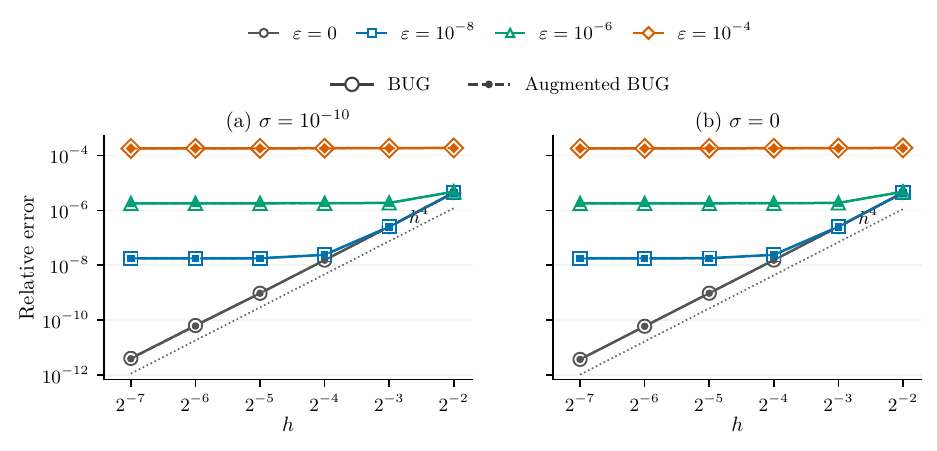}
\caption{Relative RK4-Merson errors at $T=1$ for
\eqref{s5:eq:forced-field}, with working rank $r=4$ and different
normal-component bounds $\varepsilon$. Large open markers and solid lines
denote BUG; small filled markers and dashed lines denote augmented BUG.
Dotted lines indicate fourth order.}
\label{s5:fig:normal-defect}
\end{figure}

To examine the contribution of a nonzero normal component, we also consider
\begin{equation}\label{s5:eq:forced-field}
 F_\varepsilon(t,Y)=F(t,Y)+\varepsilon C,
 \qquad \|C\|_2=1,
\end{equation}
where $F$ is given by \eqref{s4:eq:numerical-field} and $C$ is a fixed
matrix. This field satisfies \eqref{s4:eq:normal} with
$\varepsilon_r=\varepsilon$. Figure~\ref{s5:fig:normal-defect} shows the
RK4-Merson errors for small-singular-value and rank-deficient initial data.
Fourth-order decay is observed while the time-discretization error dominates.
Further refinement leads to an error plateau whose size decreases with
$\varepsilon$, consistently with the bound $O(h^4+\varepsilon_r)$.
BUG and augmented BUG have similar errors in both regimes.

\FloatBarrier

\subsection{Kinetic shear Alfv\'en waves}
\label{s5:sec:alfven}
We next study the performance of BUG without augmentation for kinetic shear
Alfv\'en waves \cite{einkemmer2023alfven}, comparing its wave dynamics and
conservation errors with augmented BUG. The distribution function satisfies
\begin{equation}\label{s5:eq:alfven}
 \partial_t f+v\partial_z f+
 \frac{\partial_z\phi+\partial_t A_\parallel}{M_e}\partial_v f=0,
\end{equation}
where the scalar and vector potentials are determined self-consistently.
We use a Fourier discretization on a $4\times4\times8\times512$ grid.
The partition $(x,y)\mid(z,v)$ gives a $16\times4096$ matrix, with initial value
\[
 f(0,x,y,z,v)=\sqrt{\frac{M_e}{\pi}}e^{-M_ev^2}
 \bigl(1+\alpha\cos(k_xx)\cos(k_yy)\cos(k_zz)\bigr).
\]
We use the weak-wave parameters of \cite{einkemmer2023alfven}, with
perturbation amplitude $\alpha=10^{-3}$. The initial matrix has rank two.
All four methods use the same spatial discretization, working rank $r=5$,
and step size $h=1.25\times10^{-5}$, and integrate to $T=0.1$.
The time discretizations are forward Euler and RK4-Merson.

The transport and acceleration terms are evaluated in factored form.
Their factors are stored at each stage and reused for the core projections
as described in Section~\ref{s5:sec:cost}.

\begin{figure}[htbp]
\centering
\begin{subfigure}{0.48\textwidth}
 \includegraphics[width=\linewidth]{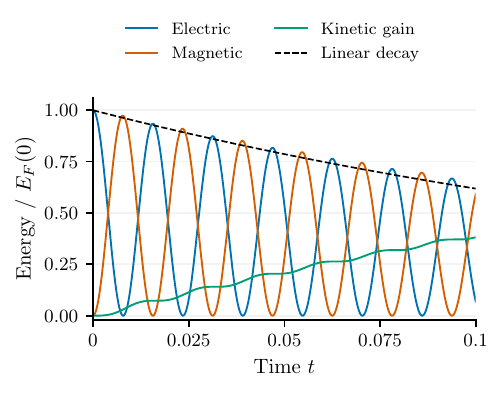}
 \caption{Euler, BUG.}
\end{subfigure}\hfill
\begin{subfigure}{0.48\textwidth}
 \includegraphics[width=\linewidth]{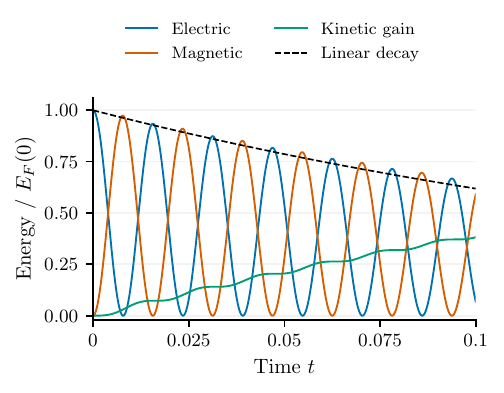}
 \caption{Euler, augmented BUG.}
\end{subfigure}

\medskip
\begin{subfigure}{0.48\textwidth}
 \includegraphics[width=\linewidth]{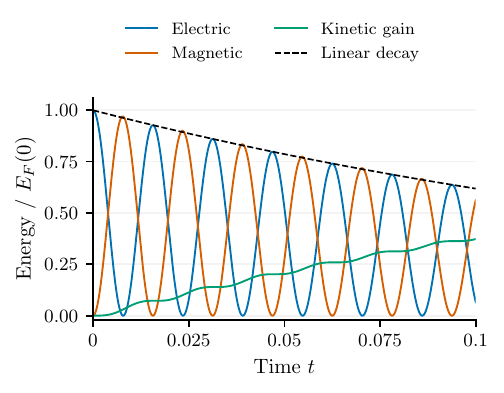}
 \caption{Merson, BUG.}
\end{subfigure}\hfill
\begin{subfigure}{0.48\textwidth}
 \includegraphics[width=\linewidth]{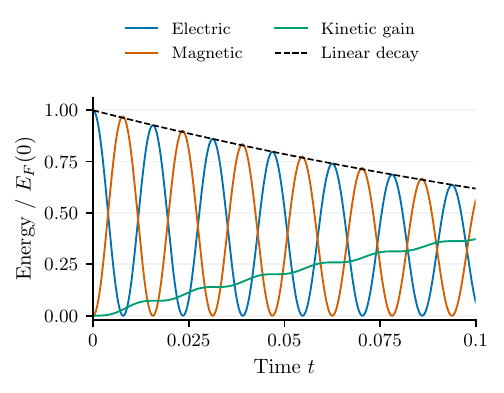}
 \caption{Merson, augmented BUG.}
\end{subfigure}
\caption{Kinetic shear Alfv\'en waves with $r=5$ and
$h=1.25\times10^{-5}$. Electric energy, magnetic energy, and kinetic-energy gain
$E_K(t)-E_K(0)$, divided by $E_F(0)$. The dashed curve is $e^{-4.8032t}$.
}
\label{s5:fig:alfven}
\end{figure}

Following \cite{einkemmer2023alfven}, Figure~\ref{s5:fig:alfven} shows the
electric and magnetic field energies together with the kinetic-energy gain.
Let $E_F$ denote the sum of the two field energies, $E_K$ the kinetic energy,
and $E=E_F+E_K$ the total energy. We normalize the energies by $E_F(0)$ to
resolve their changes on the scale of the initial wave. The dashed curve is
the linear decay envelope $e^{-4.8032t}$, using the damping rate $2.4016$
from \cite{einkemmer2023alfven}.

BUG without augmentation captures the oscillatory exchange between electric
and magnetic energies and the accompanying kinetic-energy gain. Euler
underestimates the damping, whereas Merson follows the linear decay
envelope more closely. For each time discretization, the energy histories
from BUG and augmented BUG nearly coincide.

\begin{figure}[htbp]
\centering
\includegraphics[width=0.74\textwidth]{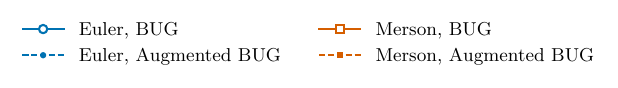}

\begin{subfigure}{0.32\textwidth}
 \includegraphics[width=\linewidth]{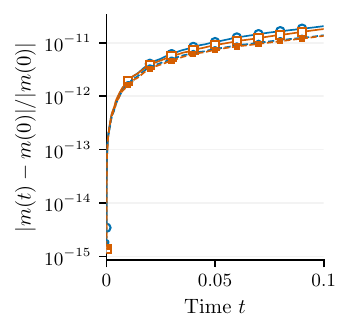}
 \caption{Mass error.}
\end{subfigure}\hfill
\begin{subfigure}{0.32\textwidth}
 \includegraphics[width=\linewidth]{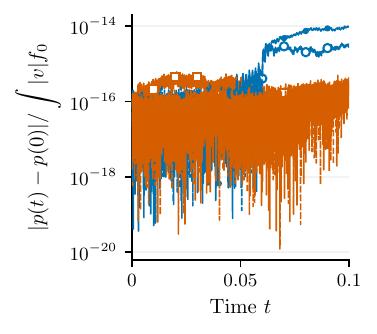}
 \caption{Momentum error.}
\end{subfigure}\hfill
\begin{subfigure}{0.32\textwidth}
 \includegraphics[width=\linewidth]{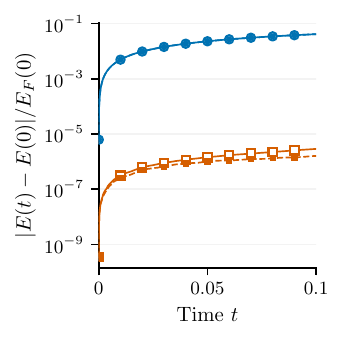}
 \caption{Total-energy error.}
\end{subfigure}
\caption{Changes in mass, momentum and total energy for the simulations in
Figure~\ref{s5:fig:alfven}, normalized by $\int f_0$, $\int |v|f_0$, and
$E_F(0)$, respectively.}
\label{s5:fig:alfven-errors}
\end{figure}

Figure~\ref{s5:fig:alfven-errors} shows small mass and momentum errors for all
four methods. The total-energy error reaches
about $4\%$ of the initial field energy with Euler, but remains at the
$10^{-6}$ level with Merson. These improvements in wave evolution and
energy balance are obtained without augmentation. The results are
comparable to those of augmented BUG in this test.

\FloatBarrier

\section{Conclusion}
\label{s4:sec:conclusion}

We have developed robust high-order BUG integrators without basis
augmentation or rank truncation. The basis and core updates use the same
prescribed increment. The two-sided relaxed-exactness estimate controls
the output error by the increment error under the column and row range
conditions. Under the assumptions on the vector field and its normal
component, increments with uniform error
$d_h=O(h^{p+1}+h\varepsilon_r)$ give the global bound
\[
 \max_{t_0\le t_k\le T}\|Y_k-A(t_k)\|
 =O(\delta+\varepsilon_r+h^p),
\]
where $\delta$ bounds the initial error and $\varepsilon_r$ bounds the
normal component. The constants and admissible step-size threshold are independent
of small singular values. The estimates also allow the actual rank of the
approximations to be smaller than $r$.

For explicit Runge--Kutta increments, the coefficient conditions in
\cite{matrix} give the stage accuracy required by this analysis. The
explicit midpoint method, Heun's third-order method, and Merson's
fourth-order method therefore have robust orders two, three, and four,
respectively. Every internal stage and final approximation has rank at
most $r$, and neither augmentation nor truncation is needed.

The same analysis allows arbitrary basis augmentation that retains the
BUG sampling spaces. For general increments, the output error and the
best rank-$r$ SVD truncation residual are controlled by the increment
error and the normal-component bound. The resulting Runge--Kutta methods
therefore retain their robust convergence orders with truncation at the
internal stages and the endpoint, without a separate truncation-error
assumption. For the same input and increment, if the augmented
spaces contain those selected by the unaugmented step, augmentation
followed by this truncation does not increase the approximation error to
$Y_0+\Delta A$. Old basis vectors or prescribed moment directions can
therefore be included without changing the convergence order.

The numerical experiments confirm the predicted orders for small singular
values and rank-deficient initial data. In the kinetic shear Alfv\'en
test, BUG without augmentation captures the wave dynamics and gives
results close to augmented BUG. Merson's method reduces the energy error
substantially compared with forward Euler in both versions.

\bibliographystyle{siamplain}
\bibliography{references}
\end{document}